\documentclass[12pt, fleqn, a4paper]{article}
\usepackage{amsthm}
\usepackage{amsmath}
\usepackage{amsfonts}
\usepackage{amssymb}
\usepackage{xspace}
\usepackage{xr}
\usepackage{float}
\usepackage{makeidx}
\usepackage{graphicx} 
\usepackage{epsfig} 
\usepackage{url}

\newcommand{\reseteqn}{\setcounter{equation}{0}}

\newcommand{\ic}{\in\mspace{-1mu}\cdot\mspace{1mu}}

\newcommand{\defas}{\mathrel{\raise.095ex\hbox{:}\mkern-4.2mu=}}
\renewcommand{\emptyset}{\varnothing}

\newcommand{\bdm}{\begin{displaymath}}
\newcommand{\edm}{\end{displaymath}}

\newcommand{\bP}{\ensuremath{\mathbf{P}}\xspace}

\newcommand{\bR}{\ensuremath{\mathbb{R}}\xspace}

\newcommand{\cE}{\ensuremath{{\cal E}}\xspace}

\newcommand{\cD}{\ensuremath{{\cal D}}\xspace}

\newcommand{\cF}{\ensuremath{{\cal F}}\xspace}
\newcommand{\cH}{\ensuremath{{\cal H}}\xspace}

\newcommand{\ra}{\ensuremath{\rightarrow}}

\renewcommand{\leq}{\ensuremath{\leqslant}}
\renewcommand{\geq}{\ensuremath{\geqslant}}

\newcommand{\hX}{\ensuremath{\hat{X}}\xspace}

\newcommand{\hZ}{\ensuremath{\hat{Z}}\xspace}

\newcommand{\mbf}{\mathbf}

\newlength{\Litil} 
\newlength{\Mid} 
\newlength{\Stor} 
\newtheorem{theorem}{Theorem}

\theoremstyle{remark}

\theoremstyle{remark}

\newcommand{\bN}{\ensuremath{\mathbb{N}}\xspace}
\newcommand{\given}{\ensuremath{\;|\;}}
\newcommand{\asdef}{\mathrel{=\mkern-4.2mu\raise.095ex\hbox{:}}}

\newcommand{\bil}{\ensuremath{\quad}}

\newcommand{\bbil}{\ensuremath{\qquad}}

\newlength{\myownsep}
\makeatother

\usepackage[all]{xy}

\begin{document}
\reseteqn

\newcommand{\ci}[1]{\ensuremath{{#1}^{\circ}}}
\newcommand{\cik}[1]{\ensuremath{{#1}^{[\circ]}}}

\title{Convergence  in Density in Finite Time Windows 
and  the
Skorohod Representation}
\author{Hermann Thorisson,
University of Iceland 
}

\date{}

\maketitle
\begin{abstract} 

\noindent 
According to the Dudley-Wichura extension of the Skorohod representation theorem,
convergence in distribution to a limit in a separable
set
is equivalent to
the existence of a coupling with elements converging a.s.\! in the metric. 
A
density analogue of this theorem says that a sequence of probability
densities on a general measurable space has a probability density as a
pointwise lower limit if and only if there exists a coupling with elements
converging a.s.\! in the discrete metric. 
In this paper the discrete-metric theorem
is extended to 
stochastic processes considered in a
widening time window. The extension is then 
used to 
 prove the separability
 version of the 
 Skorohod representation
 theorem.

\end{abstract}

\renewcommand{\thefootnote}{{}}
\footnote{\hspace{-0.8cm}
\newline
\hspace{0.8cm}{\bf Key words
} Skorohod representation, 
convergence in distribution, convergence in density,
widening time windows.
\newline {\bf Address }
Department of Mathematcs, University of Iceland,  
Dunhaga 5, 107 Reykjavik, 
Iceland, Fax: +354 552 1331, e-mail: hermann@hi.is
}

\section{Introduction}

Let $X_1, X_2, \dots, X$ be random elements 
in a general space $(E,\cE)$
with distributions 
$P_1, P_2, \dots, P$.
Let $f_1, f_2, \dots, f$ be the densities  of $P_1, P_2, \dots, P$
with respect to some  measure $\lambda$ on $(E,\cE)$.
Note that such a measure $\lambda$ always exists, we~could for instance
take $\lambda = P + \sum_{n=1}^{\infty} 2^{-n}P_n$.
If \begin{gather} \label{a}
\liminf_{n \to \infty} f_n = f\quad \text{a.e. } \lambda
\end{gather}
we  write
\begin{gather*}
\text{$X_n \to X$ in density
as $n \to \infty$.}
\end{gather*}
Note that $f_n/f$ is defined  almost everywhere $P$. 
It is the Radon-Nikodym derivative $\text{$d$}P_n/\text{$d$}P$ of 
the absolutely continuous part of $P_n$ with respect to~$P$.
Thus convergence in density 
does not depend on $\lambda$ and is  equivalent to
\begin{gather*} 
\liminf_{n \to \infty} \,\text{$d$}P_n/\text{$d$}P \,=\, 1\quad \text{a.e. } P.
\end{gather*}
In general, 
$\liminf_{n \to \infty} f_n = f \text{ a.e. } \lambda$
is  weaker than 
$\lim_{n \to \infty} f_n = f \text{ a.e. } \lambda$
and stronger than convergence in total variation.
However, if 
$(E, \cE)$ is discrete
(that is, if $E$ is countable and $\cE = 2^E \!=$ 
the class of all subsets of $E$)
then 
these 
three modes of convergence are equivalent 
and simplify to 
\begin{gather*}
\text{
$\lim_{n\to\infty}\bP(X_n = x) = \bP(X = x)$, 
\quad $x \in E$;}
\end{gather*}
see Theorems 6.1 and 7.1 in Chapter 1 of \cite{Tho:00}.

Let $(\hat{X}_1, \hat{X}_2, \dots, \hat{X})$ denote a coupling
 of $X_1, X_2, \dots, X$; this means that 
 the  random elements 
$\hat{X}_1, \hat{X}_2, \dots, \hat{X}$ 
 are defined on a common probability space
 and have the marginal distributions $P_1, P_2, \dots, P$.
In a 
 1995 paper~\cite{Tho:95}, Section 5.4, this author showed 
 that convergence in density is 
equivalent to 
the existence of a coupling  
converging in the discrete metric, that is,
\eqref{a} holds if and only if there exists
a coupling $(\hat{X}_1, \hat{X}_2, \dots, \hat{X})$
 of $X_1, X_2, \dots, X$
such that for some random variable $N$ taking values in $\bN = \{1,2,\dots\}$, 
\begin{gather}\label{b}
\hat{X}_n = \hat{X}, \quad n \geq N.
\end{gather}
This 
density result is analogous to the Skorohod representation theorem 
which says  that 
convergence in distribution 
on 
a complete separable metric  
$E$
with $\cE$ the Borel sets (a Polish space) is
equivalent to the existence of a coupling converging a.s.\!\!
in the metric.
Skorohod proved this theorem in the 1956 paper \cite{Sko:56},
Dudley removed the completeness assumption in the 1968
paper~\cite{Du:68}, and Wichura showed 
in the 1970 paper~\cite{Wi:70} that it is enough 
that the limit probability measure $P$ is concentrated 
on a separable Borel set;
for historical notes, see 
\cite{Du:02}. 
The density result in \cite{Tho:95} 
was rediscovered by Sethuraman \cite{Se:02} in 2002.
For recent~developments 
going beyond separability
and considering convergence in probability,
see the series of papers \cite{BPR:07}--\cite{BPR:14} by 
Berti, Pratelli and Rigo.

In the present paper we extend the above mentioned density result~at~\eqref{b}
to stochastic processes 
considered in a
widening time window. 
The main result, Theorem~1, is established in Section 2 while
Section~3 contains corollaries 
elaborating on that result.
In Section 4, we
show how this yields a 
new proof of
the separability version of the Skorohod representation theorem.

\section{Convergence in a widening time window}

In this section we
consider
continuous-time  
stochastic processes without restriction
on state space or paths. 
Also we allow the state space to vary with time
and include infinity in the time set. 
Discrete-time processes are considered at the end of the section.

Let $(E^t, \cE^t)$, $t\in[0,\infty]$, be a family of measurable spaces.
Let $H$ be a non-empty subset of the product set $\{(z^s)_{s\in [0,\infty]} :z^s \in E^s, s\in [0,\infty]\}$
and let $\cH$ be the smallest $\sigma$-algebra on $H$
making the maps taking $(z^s)_{s\in [0,\infty]}\in H$ to $z^t \in E^t$
measurable for all $t \in [0,\infty]$.
For $t\in[0,\infty)$, let $(H^t,\cH^t)$ be
the image space of $(H, \cH)$ under the map taking
$(z_s)_{s\in [0,\infty]} \in H$ to $
(z_s)_{s\in [0,t)}$.

If $\mbf{Z} = (Z^s)_{s \in [0, \infty]}$ 
is a random element in $(H, \cH)$
write $\mbf{Z}^t = (Z^s)_{s \in [0, t)}$ for
a segment of $\mbf{Z}$ in a finite time window of length $t \in [0, \infty)$.
Note that $\mbf{Z}^t$ is a random element in $(H^t,\cH^t)$.

According to the following theorem, convergence in density 
 in all finite time windows is the distributional form of
discrete-metric convergence in a widening time window.
(Note that the coupling in this theorem
is not a full coupling of $\mbf Z_1, \mbf Z_2, \dots, \mbf{Z}$ 
but only a coupling of $\mbf Z_1^{t_1}, \mbf Z_2^{t_2}, \dots,\mbf{Z}$.
Extensions to a full coupling are considered in the next section.)

\begin{theorem}
\label{T:3}
Let  $\mbf Z_1, \mbf Z_2, \dots, \mbf{Z}$
be  random elements in $(H, \cH)$
where $(H, \cH)$ is as above. 
Then
\begin{gather}  \label{3}
 \forall t \in [0, \infty): \bil  \mbf Z_n^t \to \mbf Z^t  \,\, \text{in density as} 
 \,\, n \to \infty\end{gather}
if and only if there exists a sequence of numbers
$0 \leq t_1 \leq t_2 \leq \dots \to \infty$ 
and a coupling 
$(\mbf{\hat{Z}}_1^{t_1}, \mbf{\hat{ Z}}_2^{t_2}, \dots, \mbf{\hat{Z}})$
of $\mbf Z_1^{t_1}, \mbf Z_2^{t_2}, \dots, \mbf{Z}$ such that 
for  some $\bN$-valued random variable $N$, 
\begin{gather}\label{r}
  \mbf {\hZ}_n^{t_n} = \mbf {\hZ}^{t_n}, \bil n \geq N.
  \qquad 
\end{gather}
\end{theorem}
\noindent
{\sc Proof.}
First, assume existence of the coupling. 
Fix $t \in [0, \infty)$, take $m \in \bN$
such that $t_{m} \geq t$, and note that then \eqref{r} yields
 $\mbf {\hZ}_n^{t} = \mbf {\hZ}^{t}$ for
$n \geq \max\{N, m\}$.
Use this and the fact that \eqref{b} implies convergence in density to obtain 
\eqref{3}.

Conversely assume that \eqref{3} holds. 
With $t \in [0, \infty)$ and $n \in \bN$,
let $Q$ be the distribution of $\mbf Z$, let  $Q^t$ be the distribution of $\mbf Z^t$,
 let $Q_n^t$ be the distribution of $\mbf{Z}_n^t$,
let $f_n^t$ be the density of $\mbf{Z}_n^t$ 
with respect to some measure $\lambda^t$
on $(H^t, \cH^t)$, and let 
$\nu_n^t$ be the measure
on $(H^t, \cH^t)$ with density $g_n^t := \inf_{k \geq n}f_k^t$. 
Due to the assumption \eqref{3}, $g_n^t$ increases
to a density of $\mbf{Z}^t$ as $n \to \infty$.
Thus by monotone convergence,
the measure $\nu_n^t$  increases setwise to $Q^t$,
\begin{gather*}
\nu_1^t \leq \nu_2^t \leq \dots \nearrow Q^t,\quad 
 t \in [0, \infty).
\end{gather*}
Thus there are numbers  $0 \leq t_1 \leq t_2 \leq \dots \to \infty$ 
such that
\begin{gather*}
0\leq Q^{t_n} - \nu_n^{t_n}  \leq 2^{-n},  \quad n\in\bN. 
\end{gather*}
For \,$A\in  \cH$ 
and $\mbf{z}^{t_n} \in H^{t_n}$, let
 \,$q^{t_n}(A\,|\,\mbf{z}^{t_n})$ be the conditional probability of the event 
\,$\{\mbf Z \in A\}$ given 
$\mbf Z^{t_n} = \mbf{z}^{t_n}$. 
Then 
\begin{gather*}
Q(A) = \int q^{t_n}(A\,|\,\cdot)\, \text{$d$} Q^{t_n},\quad 
A \in \cH.
\end{gather*}
Since $\nu_n^{t_n} \leq Q^{t_n}$ the measure
$\nu_n^{t_n} $ is absolutely continuous with respect to $Q^{t_n}$.
Thus we can
extend  $\nu_n^{t_n}$ 
from $(H^{t_n}, \cH^{t_n})$ to a measure
$\nu_n$ on  
$(H, \cH)$ by
\begin{gather*}
\nu_n(A) \defas \int q^{t_n}(A\,|\,\cdot)\, \text{$d$} \nu_n^{t_n},\quad A \in \cH.
\end{gather*}
The last three displays yield
\begin{gather*}
0 \leq Q - \nu_n \leq 2^{-n},\quad n \in \bN.
\end{gather*}
Let $h_n$ be a density of $\nu_n$ with respect to  $Q$.
For integers $n < m$ let $\mu_{n,m}$ be the measure with 
density $\min_{n \leq j \leq m}h_j$ with respect to $Q$.
Partition $H$ into sets 
$A_n, \dots , A_m \in \cH$
such that $\min_{n \leq j \leq m}h_j = h_i$ on $A_i$ and thus
\begin{gather*}
\mu_{n,m}(\cdot \cap A_i) = \nu_i(\cdot \cap A_i), \quad n \leq i \leq m.
\end{gather*}
The last two displays yield
\begin{gather*}
0 \leq Q - \mu_{n,m} = \sum_{i=n}^{m} (Q(\cdot \cap A_i) - \nu_i(\cdot \cap A_i)) \leq  
\sum_{i=n}^\infty 2^{-i} = 2^{-n+1}.
\end{gather*}
Let $\mu_n$ be the measure with density 
$\inf_{n \leq i < \infty}h_i$ with respect to $Q$ and 
send $m \to \infty$ to obtain $0 \leq Q - \mu_n  \leq  2^{-n+1}$.
Thus  the $\mu_n$ increase setwise to $Q$,
\begin{gather}\label{Q}
0 =: \mu_0 \leq \mu_1 \leq \mu_2 \leq \dots \nearrow Q. 
\end{gather}
Let  \,$\mu_n^t$\, be
the marginal of $\mu_{n}$ on $(H^t, \cH^t)$.
Note that \,$\nu_n^{t_n}$\, is the marginal of~$\nu_n$ on $(H^{t_n}, \cH^{t_n})$
and that $\mu_n \leq \mu_{n,m}$ and $\mu_{n,m} \leq \nu_n$.
Thus $\mu_n^{t_n} \leq \nu_n^{t_n}$.
Now~$\nu_n^{t}$ has density $g_n^t = \inf_{k\geq n}f_k^t$
 and $Q_n^{t}$ has density $f_n^{t}$ and 
thus $\nu_n^{t} \leq Q_n^{t}$. Since
 $\mu_n^{t_n} \leq \nu_n^{t_n}$ this yields
\begin{align}\label{Qt}
\mu_n^{t_n} \leq  Q_n^{t_n}, \quad n \in \bN.
\end{align} 
Keep in mind \eqref{Q} and \eqref{Qt} 
throughout the following coupling construction.

Let $(\Omega, \cF, \bP)$ be a probability space supporting 
the following collection of independent random elements 
with distributions to be specified below:
\begin{align*}
\text{$N$, $\mbf{V}_1$, $\mbf{V}_2,\dots$,  
$\mbf{W}_1$, $\mbf{W}_2, \dots$}
\end{align*}
Let $N$ 
be $\bN$-valued with 
distribution function 
(see \eqref{Q})
\begin{align*}
\bP(N \leq n) = \mu_{n}(H),
\quad n \in \bN.
\end{align*}
Let\, $\mbf{V}_n$\,
be\, a\, random\, element\, in\, $(H, \cH)$ with~distribution (see \eqref{Q})
\begin{gather*}
\text{$\frac{\mu_n - \mu_{n-1}}{\bP(N = n)}$\qquad 
(arbitrary distribution if $\bP(N = n) = 0$).}
\end{gather*}
Let $\mbf W_n$ be a random element in $(H^{t_n}, \cH^{t_n})$ with
distribution  (see \eqref{Qt})
\begin{gather*}
\text{$\frac{Q_n^{t_n} - \mu_{n}^{t_n}}{\bP(N > n)}$\qquad 
(arbitrary distribution if $\bP(N > n) = 0$).}
\end{gather*}
Put 
$\mbf \hZ = \mbf V_N$ to obtain that $\mbf \hZ$ has the 
same distribution as $\mbf Z$,
\begin{align*}
\bP(\mbf{ \hat{ Z}} \in \cdot) 
= \sum_{n=1}^{\infty}
\bP(\mbf{V}_n \in \cdot)\bP(N = n)   
=\sum_{n=1}^{\infty} (\mu_n - \mu_{n-1}) 
= Q.
\end{align*}
Put
$\mbf{\hat{Z}}_n^{t_n} = \mbf{V}_N^{t_n}$ on $\{N \leq n\}$ and
$\mbf{\hat{Z}}_n^{t_n}  = \mbf{W}_n$ on $\{N > n\}$
to obtain that $\mbf{\hat{ Z}}_n^{t_n}$ 
has the same distribution as $\mbf Z_n^{t_n}$,
\begin{align*}
\bP(\mbf{ \hat{ Z}}_n^{t_n} \in \cdot) 
&= \sum_{k=1}^{n}
\bP(\mbf{V}_k^{t_n}\in \cdot)\bP(N = k) + \bP(\mbf{W}_n \in \cdot)\bP(N > n)
\\&=   \sum_{k=1}^{n} (\mu_k^{t_n} - \mu_{k-1}^{t_n}) + (Q_n^{t_n} - \mu_{n}^{t_n}) 
= Q_n^{t_n}.
\end{align*}
By definition $\mbf \hZ = \mbf V_N$ and thus $\mbf \hZ^{t_n} = \mbf V_N^{t_n}$.
Also by definition, $\mbf{\hat{Z}}_n^{t_n}  = \mbf V_N^{t_n}$ 
 on $\{N \leq n\}$. Thus
$\mbf{\hat{Z}}_n^{t_n}  = \mbf{\hat Z}^{t_n}$ when $n \geq N$, 
that is, \eqref{r} holds.
\qed

\vspace{5mm}

If $\mbf{Z} = (Z^1,Z^2, \dots, Z^\infty)$
write $\mbf{Z}^k = (Z^1,Z^2, \dots, Z^k)$ for
a segment in a finite time window of length $k \in \bN\cup\{0\}$.
The following is a discrete-time version of Theorem 1.

\vspace{3mm}
\noindent
{\bf Corollary 1.}
\emph{
Let $\mbf Z_1, \mbf Z_2, \dots, \mbf{Z}$
be random elements in some product space 
$(E^1,\cE^1)\otimes(E^2, \cE^2)\otimes\dots\otimes(E^\infty, \cE^\infty)$.
Then
\begin{gather*}  
 \forall k \in \bN: \bil  \mbf Z_n^k \to \mbf Z^k \,\, \text{in density as} 
 \,\, n \to \infty
 \end{gather*}
if and only if there exists a sequence of integers
$0 \leq k_1 \leq k_2 \leq \dots \to \infty$ 
and a coupling 
$(\mbf{\hat{Z}}_1^{k_1}, \mbf{\hat{ Z}}_2^{k_2}, \dots, \mbf{\hat{Z}})$
of $\mbf Z_1^{k_1},  \mbf Z_2^{k_2}, \dots, \mbf{Z}$ such that 
for some $\bN$-valued random variable $N$, 
\begin{gather*}
  \mbf {\hZ}_n^{k_n} = \mbf {\hZ}^{k_n}, \bil n \geq N.
\end{gather*}
}

\vspace{-3mm}
\noindent
{\sc Proof.} Apply Theorem 1 to $(Z_{1}^{[1+s]})_{s\in [0,\infty]}, 
(Z_{2}^{[1+s]})_{s\in [0,\infty]},\dots, (Z^{[1+s]})_{s\in [0,\infty]}$.
(Or  repeat the proof of
Theorem 1 with
$t$ and $t_n$ replaced by  $k$ and $k_n$.)
\qed

\section{Extensions to a full coupling}

The coupling in Theorem 1 
is not a full coupling of $\mbf Z_1, \mbf Z_2, \dots, \mbf{Z}$ 
but only a coupling of $\mbf Z_1^{t_1}, \mbf Z_2^{t_2}, \dots,\mbf{Z}$.
However, in the discrete-time 
case of Corollary~1, if we 
restrict all but the infinite-time state space 
to be
discrete,
then there is 
the following simple
extension of the coupling. It~will be used in Section~4 to 
establish the separability version of the Skorohod representation theorem.

\vspace{3mm}
\noindent
{\bf Corollary 2.}
\emph{\!\!
Let $\mbf Z_1, \mbf Z_2, \dots, \mbf{Z}$
be~random elements in the product space
$(E^1,\cE^1)\otimes(E^2, \cE^2)\otimes\dots\otimes(E^\infty, \cE^\infty)$
where $(E^1,\cE^1),(E^2, \cE^2), \dots$ are discrete
and $(E^\infty, \cE^\infty)$ is some measurable space.
Then
\begin{gather*}  
 \forall k \in \bN: \bil  \mbf Z_n^k \to \mbf Z^k \,\, \text{in density as} 
 \,\, n \to \infty\end{gather*}
if and only if there exists 
a coupling 
$(\mbf{\hat{Z}}_1, \mbf{\hat{ Z}}_2, \dots, \mbf{\hat{Z}})$
of \,$\mbf Z_1,  \mbf Z_2, \dots, \mbf{Z}$ 
such~that,
for some 
$\bN$-valued random variable $N$\! 
and integers $0 \leq k_1 \leq k_2 \leq \dots \to \infty$,
\begin{gather*}
  \mbf {\hZ}_n^{k_n} = \mbf {\hZ}^{k_n}, \bil n \geq N.
\end{gather*}
}

\vspace{-3mm}
\noindent
{\sc Proof.} 
Due to Corollary 1, we only need to show that
 $(\mbf{\hat{Z}}_1^{k_1}, \mbf{\hat{ Z}}_2^{k_2}, \dots, \mbf{\hat{Z}})$
can be extended to a coupling of $\mbf Z_1,  \mbf Z_2, \dots, \mbf{Z}$.
For that purpose set, for  $n\in \bN$ and $\mbf{i}^{k_n} \in 
E^{1} \times E^{2} \times \dots \times E^{k_n}$,

\vspace{-7mm}

\begin{gather}\label{cond}
\text{$Q_{n,\mbf{i}^{k_n}} =$ the conditional distribution of $\mbf{Z}_n$ 
given 
$\{\mbf{Z}_n^{k_n}= \mbf{i}^{k_n}\}$.}
\end{gather}

\vspace{-2 mm}

\noindent
Let the probability space
$(\Omega, \cF, \bP)$  supporting 
 $\mbf{\hat{ Z}}_1^{k_1},\mbf{\hat{ Z}}_2^{k_2}, \dots$, $\mbf{\hat{Z}}$, $N$
  be large enough to also support 
  random elements in 
 $(E^1,\cE^1)\otimes(E^2, \cE^2)\otimes\dots\otimes(E^\infty, \cE^\infty)$,
\begin{gather*}
 \mbf {V}_{n,\mbf{i}^{k_n}} , \quad   n \in \bN,\,\, \mbf{i}^{k_n} \in 
 E^1\times E^2\times\dots\times E^{k_n},
 \end{gather*}

\vspace{-2 mm}

\noindent
that are independent,  independent of 
$(\mbf{\hat{ Z}}_1^{k_1},\mbf{\hat{ Z}}_2^{k_2}, \dots, \mbf{\hat{Z}},N)$, and 
such that 
\begin{align*}
\text{$\mbf {V}_{n,\mbf{i}^{k_n}}^{k_n}=\mbf{i}^{k_n}$ and
$\mbf {V}_{n,\mbf{i}^{k_n}}$ has distribution $Q_{n,\mbf{i}^{k_n}}$.}
 \end{align*}

 \vspace{-1 mm}

\noindent
 Define $\mbf {\hZ}_{n}=\mbf {V}_{n,\mbf{\hZ}_n^{k_n}}$. Then 
 \begin{gather*}
 \bP(\mbf{\hZ}_{n} \ic) =
 \sum_{\mbf{i}^{k_n}
}
 \bP(\mbf{V}_{n,\mbf{i}^{k_n}} \ic)\bP(\mbf{\hZ}_{n}^{k_n} = \mbf{i}^{k_n})
 =
 \sum_{\mbf{i}^{k_n}}
Q_{n,\mbf{i}^{k_n}}\bP(\mbf{\hZ}_{n}^{k_n} = \mbf{i}^{k_n}).
 \end{gather*}

\vspace{-3 mm}

\noindent
Since $\mbf{\hZ}_{n}^{k_n}$ has the same 
distribution as $\mbf{Z}_{n}^{k_n}$   we obtain from this and
\eqref{cond} that $\mbf{\hZ}_{n}$ has the same 
distribution as $\mbf{Z}_{n}$, as desired.
\qed

\vspace{3 mm}

In Corollary 2 we obtained a full coupling of 
$\mbf Z_1$,  $\mbf Z_2, \dots$, $\mbf{Z}$
in the discrete-time case by restricting  
$Z_n^k$ and $Z^k$ 
to a discrete state space for $k\in\bN$
 but without restricting the state 
space of $Z_n^{\infty}$ and $Z^{\infty}$. We shall
now much weaken 
this restriction 
at the expense 
of putting a restriction on $Z^{\infty}_n$ and $Z^{\infty}$.

\vspace{3mm}
\noindent
{\bf Corollary 3.}
\emph{Let  $\mbf Z_1, \mbf Z_2, \dots, \mbf{Z}$
be  random elements in the product of Polish spaces
$(E^1,\cE^1)\otimes(E^2, \cE^2)\otimes\dots\otimes(E^\infty, \cE^\infty)$.
 Then the coupling 
$(\mbf{\hat{Z}}_1^{k_1}, \mbf{\hat{ Z}}_2^{k_2}, \dots, \mbf{\hat{Z}})$ 
in \emph{Corollary~1} can be extended to a coupling 
$\mbf{\hat{Z}}_1, \mbf{\hat{ Z}}_2, \dots, \mbf{\hat{ Z}}$
of $\mbf Z_1$,  $\mbf Z_2, \dots$,~$\mbf{Z}$.}

\vspace{3mm}
\noindent
{\sc Proof.} Let 
$(\Omega, \cF, \bP)$ be the probability space
 supporting the random elements 
 $\mbf{\hat{Z}}_1^{k_1}$,  $\mbf{\hat{ Z}}_2^{k_2}, \dots$, $\mbf{\hat{Z}}$, $N$
 in Corollary 1.
Since a countable product of Polish spaces is Polish,
 there exist probability kernels $Q_{n}(\,\cdot\! \given \!\cdot\,)$, $n\in\bN$,
such that $Q_{n}(A \!\given\! \mbf{z}^{k_n})$ is the conditional probability 
of $\{\mbf{Z}_n \in A\}$ given 
 $\mbf{ Z}_n^{k_n}=\mbf{z}^{k_n}$,
 $A\in \cE^1\otimes\cE^2\otimes\dots\otimes\cE^\infty$
 and $\mbf{z}^{k_n}\in E^1\times E^2\times\dots\times E^{k_n}$.
According to the Ionescu-Tulcea extension theorem, 
this implies that the probability space $(\Omega, \cF, \bP)$
can be extended to support 
$\mbf{\hat{Z}_1}$,  $\mbf{\hat{ Z}_2}, \dots$ 
such that $\mbf{ \hZ}_n^{k_n}$ coincides with $\mbf{ \hZ}_n$
in the time window of length $k_n$ and 
the conditional probability 
of $\{\mbf{\hZ}_n\in A\}$ given 
 $\mbf{ \hZ}_n^{k_n}=\mbf{z}^{k_n}$ is
 $Q_{n}(A \!\given\! \mbf{z}^{k_n})$,
 $A\in \cE^1\otimes\cE^2\otimes\dots\otimes\cE^\infty$
 and $\mbf{z}^{k_n}\in E^1\times E^2\times\dots\times E^{k_n}$.
From this it follows that 
 $\mbf{ \hZ}_n$ has the same distribution as $\mbf{ Z}_n$, as desired. 
 \qed

\vspace{3mm}
The final corollary extends Corollary 3 to continuous time.

\vspace{3mm}
\noindent
{\bf Corollary 4.}
\emph{Let  $\mbf Z_1$,  $\mbf Z_2, \dots$,~$\mbf{Z}$
be random elements in 
$(D,\cD)\otimes(E, \cE)$ where 
$(D,\cD)=(D[0,\infty),\cD[0,\infty))$ is the Skorohod space 
of a Polish space and $(E, \cE)$
is Polish.
 Then the coupling 
$(\mbf{\hat{Z}}_1^{t_1}, \mbf{\hat{ Z}}_2^{t_2}, \dots, \mbf{\hat{Z}})$ 
in \emph{Theorem~1} can be extended to a coupling 
$(\mbf{\hat{Z}_1}, \mbf{\hat{ Z}_2}, \dots, \mbf{\hat{Z}})$
of $\mbf Z_1, \mbf Z_2, \dots, \mbf{Z}$.}

\vspace{3mm}
\noindent
{\sc Proof.} 
The Skorohod space $(D,\cD)$ is Polish 
 and thus the product 
 $(H,\cH)= 
 (D,\cD)\otimes(E,\cE)$ is Polish.
 Proceed as in the proof of Corollary 3 
 referring to Theorem 1 rather than Corollary 1,
 replacing $k_n$ by $t_n$, and
with $A\in 
 \cH$ and $\mbf{z}^{t_n}\in 
 H^{t_n}$.
\qed

\section{The Skorohod Representation}

In this section let $E$ be a metric space with metric $d$ and $\cE$ its Borel subsets.
Recall that $X_n$ is said to converge to $X$ in distribution as $n \to \infty$ 
if for all bounded 
continuous functions $h$ from $E$ to $\bR$,
\begin{gather*}
\int h \,\text{$d$}P_n \to \int h \,\text{$d$}P, \quad n \to \infty.
\end{gather*}
Recall also that $A \in \cE$ is 
called a $P$-continuity set 
if $P(\partial A) = 0$ where $\partial A$ denotes the boundary of $A$, and that 
by the Portmanteau Theorem (Theorem 11.1.1 in \cite{Du:02}) 
convergence in distribution is equivalent to 
\begin{gather}\label{Port}
\text{$P_n (A) \to P(A)$ as $n \to \infty$ for all $P$-continuity 
sets $A$.}
\end{gather}
We shall now use Corollary 2 to prove the Skorohod representation theorem
in the separable case.

\begin{theorem}
\label{T:1}
Let $X_1$, $X_2, \dots$, $X$ be random elements 
in 
a metric space $E$  
equipped with its Borel subsets $\cE$. 
Further, let $X$   
take values 
almost surely
in a separable subset 
$E_0 \in
\cE$. 
Then
\begin{gather}  \label{1}
  X_n \ra X \,\, \text{in distribution 
  as}
   \,\, n \to \infty
\end{gather}
if and only if there is a coupling 
$(\hX_1, \hX_2, \dots, \hX)$
of $X_1$, $X_2, \dots$, $X$ such that
\begin{gather} \label{pontwise}
 \hX_n \to \hX   \,\, \text{pointwise  as}
   \,\, n \to \infty.
\end{gather}
\end{theorem}
\vspace{3mm}
\noindent
{\sc Proof.} 
Let $d$ be the metric.
We begin with basic 
preliminaries.
First, assume existence of the coupling and let $h$ be a bounded continuous function.
Then 
(10) yields that
$h(\hX_n) \to h(\hX)$  pointwise
as $n \to \infty$ and by bounded convergence,
$\int h \,\text{$d$}P_n \to \int h \,\text{$d$}P$  as $n \to \infty$. 
Thus 
\eqref{1} holds.

Conversely, assume from now on
that \eqref{1}, and thus \eqref{Port}, holds.
For each $\epsilon > 0$, 
any separable Borel set
can be covered by countably 
many $E$-balls of diameter $< \epsilon$.
Note  that for every 
$y \in E$ and $r > 0$,
$ \partial\{x \in E : d(y,x) < r \} \subseteq {\partial\{x \in E : d(y,x) = r \}}$
and that the set on the right-hand side
has $P$-mass~$0$ except for countably many
radii~$r$.
Thus the covering sets below
may be taken to be $P$-continuity~sets.
Moreover, since 
$\partial(A \cap B) \subseteq \partial A \cup
\partial B$ for all subsets $A$ and $B$ of $E$,
the covering sets can be taken to~be disjoint. 

Let $A_2, A_3, \dots$ be 
disjoint $P$-continuity sets of diameter $< 1$ covering 
$E_0$
and put \,$A_1 = E \setminus (A_2 \cup A_3 \cup \dots$).\,
Then $A_1$ is  also a $P$-continuity set since 
$P(A_1) 
 = 0$ and since \,$\partial{A_1} $ 
cannot contain interior points
of the $P$-continuity sets $A_2, A_3, \dots$
Thus $\{A_i\!:\!i\in \bN\}$ 
 is a partition 
of  $E$ into $P$-continuity~sets. 
Put\, $A_{11}=A_1$\,
and \,$A_{12} = A_{13} = \dots = \emptyset $.\, 
For 
$i > 1$, let  $A_{i2}, A_{i3}, \dots$
be disjoint $P$-continuity subsets  of $A_i$ of diameter $< 1/2$ covering 
$E_0\cap A_i$
and put $A_{i1} = A_i \setminus (A_{i2} \cup A_{i3} \cup \dots$).\,
Then again $\{A_{\mbf{i}^2}\! : \mbf{i}^2 \in \mbf{\bN}^2\}$
is a partition 
of  $E$ into $P$-continuity sets.
Continue this recursively in $k\in \bN$ to obtain 
a sequence of 
partitions $\{A_{\mbf{i}^k}\! : \mbf{i}^k\! \in \mbf{\bN}^k\}$ 
of $E$ into $P$-continuity sets 
such that 
\begin{align}\label{E0}
 \text{
$A_{\mbf{i}^k}$, $\mbf{i}^k\! \in (\bN\setminus \{1\})^k$,
cover $E_0$ and are each
of diameter $< 1/k$}
\end{align}
and such that 
the 
partitions are nested in the sense that for $k \in \bN$ and $\mbf{i}^k\! \in \bN^k$
it holds that
$A_{\mbf{i}^{k}}=A_{\mbf{i}^{k}1} \cup A_{\mbf{i}^{k}2}\cup\dots$

After these basic 
preliminaries, we are now ready to apply Corollary~2.
Let $\mbf Z_1$,  $\mbf Z_2, \dots$, $\mbf{Z}$ be the
random elements in
$(\bN, 2^{\bN})^{\bN}\otimes (E,\cE)$ defined as follows
(well-defined because the partitions are nested):
set $Z_n^{\infty} = X_n$ and $Z^{\infty} = X$ and for $k\in\bN$
\begin{align*}
 \mbf{Z}_n^k = \mbf{i}^k \bil \text{if} \bil X_n \in A_{\mbf{i}^k} \bbil &\text{and}\bbil 
    \mbf{Z}^k =   \mbf{i}^k \bil \text{if} \bil \,\,X \in A_{\mbf{i}^k} .
\end{align*}
Due to \eqref{Port}, we have
$\bP(\mbf{Z}_n^k = \mbf{i}^k) \to \bP(\mbf{Z}^k =
\mbf{i}^k)$ as $ n  \to  \infty$, $\mbf{i}^k \in \bN^k$, $k\in \bN$.
Thus  $\mbf{Z}_n^k \to \mbf{Z}^k$ in density as $n\to \infty$
and Corollary~2 yields the existence of a coupling 
$(\mbf{\hat{ Z}}_1,\mbf{\hat{ Z}}_2, \dots$, $\mbf{\hat{Z}})$ 
of $(\mbf{Z}_1,\mbf{Z}_2, \dots$, $\mbf{Z})$,
an $\bN$-valued random variable~$N$\!, 
and integers $0 \leq k_1 \leq k_2 \leq \dots \to \infty$,
such that 
\begin{gather} \label{ZZ}
  \mbf {\hZ}_n^{k_n} = \mbf {\hZ}^{k_n}, \bil n \geq N.
\end{gather} 
 Now define the coupling of $X_1, X_2, \dots, X$
 by setting $\hX_n = \hZ_n^{\infty}$ and $\hX = \hZ^{\infty}$. 
 Then (after deleting a null event)
  we have that $\hX \in E_0$ and that
  for $k\in\bN$
\begin{gather*}
\mbf{\hZ}_n^k = \mbf{i}^k \bil \text{if} \bil \hX_n \in A_{\mbf{i}^k}\qquad\text{and}
\qquad \mbf{\hZ}^k = \mbf{i}^k \bil \text{if} \bil \hX \in A_{\mbf{i}^k}.
\end{gather*}
Thus $\hX_n \in A_{\mbf{\hZ}_n^{k_n}}$ and $\hX \in A_{\mbf{\hZ}^{k_n}}$ 
for all $n\in\bN$.
Apply \eqref{ZZ} to obtain that
\begin{gather}\label{AA}
\text {both \,$\hX_n\in A_{\mbf{\hZ}^{k_n}}$\, and\, $\hX\in A_{\mbf{\hZ}^{k_n}}$ \,when \,$n\geq N$.}
\end{gather}
Finally, apply \eqref{E0}: since $\hX \in E_0$ we have that 
$\mbf{\hZ}^{k_n}
  \in (\mbf{\bN}\setminus \{1\})^{k_n}$
so $A_{\mbf{\hZ}^{k_n}}$
 has diameter $< 1/k_n$.
From this and \eqref{AA} we obtain that 
\begin{gather*}
d(\hX_n, \hX) < 1/k_n,  \quad 
n\geq N.
\end{gather*}
Since $N < \infty$
and $\lim_{n\to \infty} 1/k_n = 0$ 
this 
implies that  $d(\hX_n, \hX) \to 0$ pointwise, that is,
(10) holds. 
\qed

\end{document}